\documentclass[12pt,a4paper]{amsart}

\usepackage {amsfonts}
\usepackage[english]{babel}
\def\b{\beta}
\def\g{\gamma}
\def\d{\delta}
\def\a{\alpha}
\def\p{\varphi}
\def\e{\varepsilon}
\def\l{\lambda}
\def\t{\theta}
\def\R{{\mathbb R}}
\def\C{{\mathbb C}}

\def\Z{{\mathbb Z}}
\def\hn{\frak n}
\def\hR{\frak R}

\def\mod{\mbox{mod}}
\def\supp{\mbox{supp}}

\def\o{\omega}
\def\p{\varphi}
\def\L{{\frak J}}

\newtheorem{Th}{Theorem}
\newtheorem{ThC}{Theorem C}
\newtheorem*{ThL}{Theorem LO}
\newtheorem*{ThH}{Theorem H}

\newtheorem*{de}{Definition}
\newtheorem*{Cor}{Corollary}
\newtheorem{Pro}{Proposition}

\begin{document}

\title{Fourier quasicrystals  and Lagarias' conjecture}

\author{S.Yu. Favorov}

\address{Sergey Yuryevich Favorov,
\newline\hphantom{iii} The Karazin Kharkiv National University
\newline\hphantom{iii} Svoboda squ., 4,
\newline\hphantom{iii} 61022, Kharkiv, Ukraine}
\email{sfavorov@gmail.com}

\maketitle {\small
\begin{quote}
\noindent{\bf Abstract.}
J.C.Lagarias (2000) conjectured that if $\mu$ is a complex measure on p-dimensional Euclidean space with a uniformly discrete support and its spectrum (Fourier transform) is also a measure with a uniformly discrete support,
then the support of $\mu$ is a subset of a finite union of shifts of some full-rank lattice. The conjecture was proved by N.Lev and A.Olevski (2013) in the case p=1. In the case of an arbitrary p they proved the conjecture only for a positive measure $\mu$.

Here we show that Lagarias' conjecture is false in  general case and find two new special cases when assertion of the conjecture is valid.
\medskip


\end{quote}
}

\bigskip

Let $\mu$ be a  complex-valued measure in $\R^p$.  Suppose that $\mu$ is slowly increasing, i.e., its variation $|\mu|$ satisfies
the condition $|\mu|\{|x|<R\}=O(|R|^M)$ as $R\to\infty$ for some $M<\infty$. Hence, $\mu$ is a continuous linear functional in the space $\L$ of rapidly decreasing $C^\infty$-functions with seminorms
 $$
 p_N(f)=\sup_{\max_j\a_j<N}\sup_{x\in\R^p}(1+|x|^N)|D_\a f(x)|,
 $$
where $D_\a$ are partial derivatives of the order $\a_1,\dots,\a_p$.
The Fourier transform of the measure $\mu$ is defined by the equality
$\hat\mu(f)=\mu(\hat f)$ for all $f\in\L$; here
 $$
   \hat f(y)=\int_{\R^p}f(x)\exp\{-2\pi i<x,y>\}dx
 $$
is a Fourier transform of the function $f$. We will consider the case of uniformly discrete $\supp\mu$, which means
$|x-x'|\ge\g$  for all $x,x'\in\supp\mu$ and some $\g>0$. Following \cite{LO}, we will say that $\supp\mu$ is a Fourier quasicrystal, if $\supp\hat\mu$ is a pure point measure, or equivalently, if $\supp\hat\mu$ is countable (possibly dense in $\R^p$). We will say also  that $\supp\hat\mu$ is a spectrum of the quasicrystal. These notions were inspired by experimental discovery in the middle of 80's of non-periodic atomic structures with diffraction patterns consisting of sports.

Lagarias' conjecture takes its origin in the classical Poisson summation formula. Let $f$ be a sufficiently smooth and rapidly decreasing function on $\R^p$. Then
 $$
  \sum_{n\in\Z^p}f(n)=\sum_{n\in\Z^p}\hat f(n).
 $$
  In other words, the measure $\mu=\sum_{n\in\Z^p}\d_n$, where $\d_a$ means the usual Dirac measure (unit mass) at the point $a\in\R^p$, satisfies the condition $\hat\mu=\mu$. It is easy to see that for a full-rank lattice $L=A(\Z^p)$, where $A$ is
a non-degenerate linear operator in $\R^p$, and the conjugate  lattice $L^*=\{y\in\R^p: <x,y>\in\Z\}$ we get
 $$
     \widehat{\left( \sum_{x\in L}\d_x\right)}= (\det A)^{-1}\sum_{y\in L^*}\d_y.
 $$
 The converse is also true:
 \begin{ThC}[A.Cordoba \cite{C1}]
 Let $\{x_n\},\ \{y_n\}$ be uniformly discrete sets in $\R^p$, $c_n>0$ for all $n$,
  $$
  \mu=\sum_n \d_{x_n},\quad \hat\mu=\sum_n c_n\d_{y_n}.
  $$
Then there is a full-rank lattice $L=A(\Z^p)$ such that $\{x_n\}=L$, $\{y_n\}=L^*$, $c_n=(\det A)^{-1}$.
\end{ThC}

\bigskip
J.C.Lagarias (\cite{La}, p.79) conjectured that if
 $\mu$ is a complex measure on $\R^p$ with the uniformly discrete support,
and if its spectrum $\hat\mu$ also is a measure with the uniformly discrete support,
then there is a full-rank lattice $L$ and $a_1,\dots,a_N,b_1,\dots,b_{N'}\in\R^p$ such that
 $$
  \supp\mu\subset\cup_{j=1}^N (L+a_j),\qquad \supp\hat\mu\subset \cup_{j=1}^{N'} (L^*+b_j).
 $$
  In other words, the quasicrystal is a subset of a finite union of shifts of a full-rank lattice.

The most strong result in this direction was obtained by N.Lev and A.Olevskii:
 \begin{ThL}\cite{LO}
 The Lagarias' conjecture is valid in the case $p=1$, i.e., for measures on the real axis, and  in the case of
 an arbitrary $p$ and a positive  measure $\mu$ (or $\hat\mu$). Moreover, if $ \supp\mu$ satisfies the conclusion
 of the conjecture in the case $p\ge1$, then $\mu$ is of the form
 $$
 \mu=\sum_{j=1}^{N}\sum_{x\in L+a_j}P_j(x)\d_x,
 $$
 where $P_j(x)$ are finite linear combinations of exponents $e^{2\pi i<\o,x>}$.
\end{ThL}
\medskip

We prove that Lagarias' conjecture fails in general case.
Let $L=\{(\sqrt{2}m_1,m_2)\in\R^2:\ (m_1,m_2)\in\Z^2\}.$
Then  $L^*=\{(k_1/\sqrt{2},k_2)\in\R^2:\ (k_1,k_2)\in\Z^2\}.$
Recall that for $\mu_a(E)=\mu(E+a)$ we have
$$
\widehat{(\mu_a)}(y)=e^{2\pi i<a,y>}\hat\mu(y),\qquad \widehat{(e^{2\pi i<a,x>}\mu)}(y)=\hat\mu_{-a}(y).
$$
Let
$$
 \nu=\sum_{(n_1,n_2)\in\Z^2}\d_{n_1,n_2}+\sum_{(m_1,m_2)\in\Z^2}e^{m_2\pi i}\d_{\sqrt{2}m_1,m_2+1/2}
$$
Then
$$
 \hat\nu=\sum_{(n_1,n_2)\in\Z^2}\d_{n_1,n_2}+\sum_{(k_1,k_2)\in\Z^2}\frac{e^{k_2\pi i}}{\sqrt{2}}\d_{k_1/\sqrt{2},k_2-1/2}.
$$
Then
 $$
 \supp\nu=\Z^2\cup(L+(0,1/2)),\qquad\supp\hat\nu=\Z^2\cup(L^*-(0,1/2)).
 $$
Note that $\nu$ and $\hat\nu$ are real measures with masses $\pm1$ and their supports are uniformly discrete.
Furthermore, let the set $\Z^2\cup(L+(0,1/2))$ be a subset of the union of a finite number of shifts of some lattice $K$.
Then both projections of the set on the directions of the generating vectors of $K$ are uniformly discrete sets. Clearly, one of the directions is $x_1=0$. Assume that another one is $l=(\cos\t,\sin\t),\ \t\neq\pi/2$. The projection of the set on $l$ equals
 \begin{equation}\label{1}
 \{n_1\cos\t+n_2\sin\t+m_1\sqrt{2}\cos\t+m_2\sin\t+(1/2)\sin\t:\,n_1,n_2,m_1,m_2\in\Z\}
 \end{equation}
By Kronecker's theorem, the system of inequalities
 \begin{gather*}
 |t\sqrt{2}+(1/2)\tan\t|<\e (\mod\,\Z)\\  |t|<\e(\mod\,\Z)
  \end{gather*}
has an arbitrary large solution for any $\e>0$. Therefore, for any $\e>0$ there are arbitrary large integers $s,\,r$ such that
 $$
  |s\sqrt{2}+(1/2)\tan\t+r|<\e.
 $$
 In (\ref{1}) we let $n_1=-r,\,m_1=s,\,n_2=-m_2=j$ with an arbitrary integer $j$. We get the contradiction with our choice of $l$.
\bigskip

But there are some results showing  that Fourier quasicrystal  may be a finite union of shifts of {\it several} full-rank lattices.
We need the following definition.

\begin{de}A (complex) measure $\mu$ on $\R^p$ is {\it translation bounded}, if
$$
 \sup_{x_0\in\R^p}|\mu|(B(x_0,1))<\infty.
$$
\end{de}
As usually, $B(a,r)$ is a ball of radius $r$ with center at $a$, $|\mu|$ is a variation of the measure $\mu$.

\begin{ThC}[A.Cordoba \cite{C2}]\label{T3}
 Let a uniformly discrete set $\Lambda\subset\R^p$ be given as a disjoint union of $N$ subsets $\Lambda_j$,  and $\mu=\sum_{j=1}^N\sum_{x\in\Lambda_j}a_j\d_{x_n}$, where $a_j,\ j=1,\dots,N$ are complex numbers. If Fourier transform
$\hat\mu$  is a translation bounded measure with a countable support, then $\supp\mu$ is a finite union of shifts of several  full-rank lattices.
\end{ThC}

 The principal point of the proof of Cordoba's theorem is the following assertion
 \begin{Pro}\label{P1}
 Under conditions of theorem C2 there is a measure $\hn$ on the Bohr compactification  $\hR$ of $\R^p$ such that its Fourier transform $\hat{\bar\hn}$ with respect to the dual pair  $(\hR,\,\R^p)$ is a discrete measure, $\hat{\hn}(x)=1$
for $x\in\supp\mu$, and $\hat{\hn}(x)=0$ for $x\not\in\supp\mu$.
 \end{Pro}
Note that deriving  Theorem C\ref{T3} from this proposition is based on the Helson-Cohen characterisation of idempotent measures on locally compact abelian groups (\cite{R1}, Ch.3):
\begin{ThH} Let $X$ be a locally compact abelian group, $\Gamma$ be its dual group, i.e., the group of continuous characters on $X$, and $\nu$ be an idempotent (with respect to convolution) measure on $X$.  Then $\supp\hat\nu$ belongs to the the smallest ring of subsets of $\Gamma$, which contains all open cosets in $\Gamma$.
\end{ThH}

Here we prove the following stronger version of Cordoba's theorem.

\begin{Th}\label{T4}
 Let $\{x_n\}$ be a uniformly discrete set in $\R^p$, $\mu=\sum_n \mu(x_n) \d_{x_n}$, let the set
 $\{|\mu(x_n)|\}=\{\b_1,\dots,\b_N\}$ be finite, and $\hat\mu$  be a translation bounded measure with a countable support.
  Then $\supp\mu$ is a finite union of shifts of several full-rank lattices.
\end{Th}
{\bf Proof}.  We have to check that  Proposition \ref{P1} is valid under assumptions  of Theorem \ref{T4} too.

Let $\l,\rho$ be measures on $\R^p$ such that $\l(E)=\hat\mu(-E),\,\rho(E)=\overline{\hat\mu(E)}$. Hence,
Fourier transforms of  $\l$ and $\rho$ are the measure $\mu$ and the complex conjugate to $\mu$ respectively. Clearly, the measures
$\l$ and $\rho$ are translation bounded measures with countable supports. Let $\p$ be an infinitely differentiable function such that $\supp\p\subset B(0,1)$ and $\hat\p(0)=1$. Clearly, $\hat\p(x)\to 0$ as $|x|\to\infty$. Put $\l_M=M^{-p}\p(\cdot/M)\l,\ \rho_M=M^{-p}\p(\cdot/M)\rho$. Note that Fourier transforms of these measures are infinitely differentiable functions on $\R^p$. Therefore, for any point $x\in\R^p$
 $$
   \lim_{M\to\infty}\hat\l_M(x)=\lim_{M\to\infty}\hat\p(M\cdot)*\mu(x)=\mu(x),\qquad \lim_{M\to\infty}\hat\rho_M(x)=\overline{\mu}(x).
 $$
Hence if the measure $\nu_M$ is the convolution of $k$ measures $\l_M$ and $m$ measures $\rho_M$, then
        $$
        \hat\nu_M(x)\to(\mu(x))^k(\overline{\mu(x)})^m.
        $$
The same reasoning takes place if $\nu_M$ is a linear combination of such convolutions. Therefore, if we replace
$z$ by $\l_M(x)$, $\bar z$ by $\rho_M(x)$, and multiplication by convolution in  the polynomial
   $$
           P(z,\bar z)= 1-\prod_{j=1}^N(1-z\bar z/\b_j^2),
   $$
then we obtain for $M\to\infty$
 \begin{equation}\label{2}
       P(\hat\l_M,\,\hat\rho_M)(x)\to1\ \hbox{as}\ x\in\supp\mu,\qquad P(\hat\l_M,\,\hat\rho_M)(x)\to0\ \hbox{as}\ x\not\in\supp\mu.
\end{equation}
On the other hand, since $\hat\mu$ is translation bounded, we see that the total variation of the measures $\l_M,\,\rho_M$, and $\nu_M=P(\hat\l_M,\,\hat\rho_M)$ are bounded uniformly with respect to $M$. Hence, the measures $\nu_M$ have natural extension
 to the finite measures $\hn_M$ in the Bohr compactification $\hR$ of $\R^p$, with uniformly bounded total variation. Therefore there is a subsequence $M'$ such that $\hn_{M'}\to\hn$ in the weak--star topology, and we get $\hat{\hn}_M(x)\to\hat{\hn}(x)$ for all $x\in\hR$ as $M'\to\infty$. By (\ref{2}), we obtain the conclusion of Proposition \ref{P1} in this case too.
\medskip

The assertion of Lagarias' conjecture in the original form is valid  under additional conditions on quasicrystal.

Let us recall some definitions and simple properties (see, for example, \cite{M}).

\begin{de} A continuous function $f$ on $\R^p$ is  almost periodic, if
for any  $\e>0$ the set of $\e$-almost periods of $f$
  $$
  \{\tau\in\R^p:\,\sup_{x\in\R^p}|f(x+\tau)-f(x)|<\e\}
  $$
  is a relatively dense set in $\R^p$, i.e., there is  $l=l(\e)$ such that
  any ball of radius $l$ contains an $\e$-almost period of $f$.
\end{de}
The definition is equivalent to the following one: for any $\e>0$ there is a finite exponential sum $Q(x)=\sum c_n\exp\{2\pi i<x,\o_n>\}$ such that $\sup_{x\in\R^p}|Q(x)-f(x)|<\e$.

\begin{de}
A (complex) measure $\mu$ on $\R^p$ is {\it almost periodic}, if
for any continuous function $\p$ on $\R^p$ with a compact support
 the function $\int\p(x+t)d\mu(t)$ is almost periodic in $x\in\R^p$.
\end{de}

\begin{Pro}\cite{M}
Any almost periodic measure is translation bounded.
\end{Pro}

\begin{Pro}\cite{M}
Let $\mu$ and its spectrum $\hat\mu$ be translation bounded measures.
Then $\mu$ is almost periodic iff $\hat\mu$ is a discrete measure with a countable support.
\end{Pro}

Hence it is natural to change the condition "a countable spectrum" to "almost periodic measure". Here we get the following theorem

\begin{Th}
Let $\mu_1, \mu_2$ be almost periodic discrete measures on  $\R^p$ with countable supports, and
$\inf_{x\in\R^p}|\mu_1(x)|>0$, $\inf_{x\in\R^p}|\mu_2(x)|>0$.   If the
set of differences between points of $\supp\mu_1$ and $\supp\mu_2$ is discrete, then the supports
are finite unions of shifts of a unique full-rank lattice $L$, i.e.,
  there exist $c_k^j\in\R^p,\  k=1,2,\dots,r_j,$ such that $\supp\mu_j=\cup_{k=1}^{r_j}(L+c_k^j),\ j=1,2$.
 \end{Th}
 {\bf Remark}. In the case $\mu_1=\mu_2=\sum_{x\in\Lambda}\d_x$ the condition ''$\Lambda-\Lambda$ discrete'' appeared earlier in connection with so called Meyer sets \cite{M1}. Note that the name Meyer set was assigned later
by others (see \cite{Mo}).

 {\bf Proof}.   Let $\p(x)$ be a continuous function such that $\p(x)\ge0,\ \p(0)=1$, and $\supp\p\subset B(0,1)$, sett $\p_\eta(x)=\eta^{-p}\p(x/\eta)$. The sums
  $$
  S^\eta_j(x)=\sum_{t\in\supp\mu_j}\mu_j(t)\p_\eta(x+t), \quad j=1,2,
  $$
  are almost periodic functions in $x\in\R^p$.  We prove that for any $\eta>0$ there is a relatively dense set of common $\e$-almost periods of these functions.
  Indeed, using the above alternative definition of almost periodic functions, one can prove the result for two arbitrary finite exponential sums $Q_1(x)=\sum_{n=1}^N c_n\exp\{2\pi i<x,\o_n>\}$ and $Q_2(x)=\sum_{n=1}^M b_n\exp\{2\pi i<x,\sigma_n>\}$. By Kronecker's theorem, the system of inequalities
 \begin{gather*}
 |<\tau,\o_n>|<\b (\mod\,\Z),\ n=1,\dots,N\\  |<\tau,\sigma_n>|<\b (\mod\,\Z),\ n=1,\dots,M
  \end{gather*}
has a relatively dense set of solutions for any $\b>0$. For sufficiently small $\b=\b(\e)$ it implies that the
inequalities
 $$
 \sup_{x\in\R^p}|Q_1(x+\tau)-Q_1(x)|<\e, \qquad\sup_{x\in\R^p}|Q_2(x+\tau)-Q_2(x)|<\e
  $$
  are valid for each solution $\tau$ of the system.

An evident consequence follows from the proved result : there is $R<\infty$ such that any ball of radius $R$ contains at least one point of $\supp\mu_1$ and at least one point of $\supp\mu_2$. Next, there is $r>0$ such that any ball of radius $r$ contains at most one point of $\supp\mu_1$ and at most one point of $\supp\mu_2$. Indeed, if there are sequences $x_n,\,x'_n\in\supp\mu_1,\,x_n\neq x'_n $ such that $x_n-x'_n\to0$, then one can take $y_n\in\supp\mu_2$ such that $|y_n-x_n|<R,\ |y_n-x'_n|<R+1$, hence we get infinite differences $y_n-x_n$ or $y_n-x'_n$ in the ball of radius $R+1$ that contradicts the property of $\supp\mu_1-\supp\mu_2$.

 Next, since the set of differences $\supp\mu_1-\supp\mu_1$ is discrete, we see that there is $\e>0$ such that $2\e<\min\{1;r;|(a-b)-(c-d)|\}$ whenever $a,\,c\in\supp\mu_1$, $b,\,d\in\supp\mu_2$, and $|a-b|<2R+2,\,|c-d|<2R+2,\,a-b\neq c-d$.

 Without loss of generality suppose that $|\mu_1(x)|\ge1$ for all $x\in\supp\mu_1$ and $|\mu_2(x)|\ge1$ for all
 $x\in\supp\mu_2$. If $\eta<r/2$, then for any  $x\in\R^p$  both sums $S^\eta_j(x),\,j=1,2$ contain at most one nonzero term. Let $\tau$ be a common $1/2$-almost period of these sums. If $x\in\supp\mu_1$, then $S^\eta_1(x)=1$ and $S^\eta_1(x+\tau)\neq0$, therefore for any $a\in\supp\mu_1$  there is $c\in\supp\mu_1$ such that $|a+\tau-c|<\e$. The point with this property is unique, because for another $c'\in\supp\mu_1$ we have $|a+\tau-c'|\ge |c'-c|-|a+\tau-c|>r-\e>\e$. In the same way, for any $b\in\supp\mu_2$  there is a unique $d\in\supp\mu_2$ such that $|b+\tau-d|<\e$.

  Fix $a$ and put $T=c-a$. Since $|\tau-T|<\e$, we see that for any $x\in\supp\mu_1$ and any $y\in\supp\mu_2$ there are $x'\in\supp\mu_1$ and $y'\in\supp\mu_2$ such that $|x+T-x'|<2\e,\ |y+T-y'|<2\e$. We will prove that  $T$ is a common period of $\supp\mu_1$ and $\supp\mu_2$.

Suppose that $b\in\supp\mu_2$ such that $b\neq a$ and $|a-b|<2R+1$. Then there is a point $d\in\supp\mu_2$ such that
$|b+T-d|=|(a-b)-(c-d)|<2\e$. Since $|c-d|\le|a-b|+|b+T-d|<2R+2$, we obtain $a-b=c-d$
and $d=b+T$. We repeat these arguments for all $b\in\supp\mu_2$ such that $|b-a|<2R+1$ and,
after that, for all $a'\in\supp\mu_1$ such that $|a'-b|<2R+1$, then for all $b'\in\supp\mu_2$ such that $|a'-b'|<2R+1$. After a finite or countable number of steps we obtain two sets
 $$
  A_1=\{a\in\supp\mu_1:\,a+T\in\supp\mu_1\},\qquad  A_2=\{b\in\supp\mu_2:\,b+T\in\supp\mu_2\}.
 $$
If $\supp\mu_1\setminus A_1\neq\emptyset$, then set
 $$
 R_1=\inf\{|a-a'|:\,a\in A_1,\,a'\in\supp\mu_1\setminus A_1\}.
 $$
If $R_1\ge 2R+1$, take $a\in A_1$ and $a'\in\supp\mu_1\setminus A_1$ such that $|a'-a|<R_1+1$. Then there is a point $c\in B((a+a')/2,R)\cap\supp\mu_1$. It is easy to see that $|c-a|<R_1$ and $|c-a'|<R_1$, therefore $c\not\in A_1$ and $c\not\in\supp\mu_1\setminus A_1$, which is impossible. Thus we have $R_1<2R+1$. In this case take $b\in B((a'+a)/2,R)\cap\supp\mu_2$. Since $|b-a|\le
|b-(a+a')/2|+|(a+a')/2-a|<R+R_1/2+1/2<2R+1$, we see that $b\in A_2$. On the other hand, $|b-a'|<2R+1$ as well, hence, $a'\in A_1$.
This contradiction implies that $A_1=\supp\mu_1$. In the same way, $A_2=\supp\mu_2$. Hence, $T$ is a common period of $\supp\mu_j,\,j=1,2$.

 Next, consider $p$ cones
 $$
 C_j=\left\{x\in\R^p:\,|x-\langle x,e_j\rangle e_j|<\g|x|\right\},
 \quad j=1,\dots,p,
 $$
where $e_j,\,j=1,\dots,p,$ is the intrinsic basis for $\R^p$ . There are $(1/2)$-almost periods $\tau_j\in C_j$
and, therefore, common periods $T_j\in C_j,\,j=1,\dots,p$. We may suppose that $|T_j|>1$ and $\g$ is small  enough , then
$T_j$ are linearly independent over $\R$. Consequently, the set $L=\{n_1T_1+\dots
+n_pT_p:\,n_1,\dots,n_p\in\Z\}$ is a full--rank lattice. Next, the set $F_1=\{a\in\supp\mu_1:
|a|<|T_1|+\dots+|T_p|\}$ is finite. All vectors $t\in L$ are periods of $\supp\mu_1$, hence,
$L+F_1\subset\supp\mu_1$. On the other hand, for each $a\in\supp\mu_1$ there is $t\in L$ such that
$|a-t|<|T_1|+\dots+|T_p|$, hence, $a-t\in F_1$. In the same way, there is a finite set $F_2$ such that $\supp\mu_2=L+F_2$.  The theorem is proved.
\medskip

In particular, in the case $\mu_1=-\a\mu_2$ we get the following result:
\begin{Cor}
Let $\mu$ be an almost periodic measure on  $\R^p$ with a countable support, and
$\inf_{x\in\R^p}|\mu(x)|>0$.   If the set $\{x+\a x':\,x,x'\in\supp\mu\}$ for some $\a\in\C$ is discrete, then the $\supp\mu$
is a finite union of shifts of a unique full-rank lattice $L$.
 \end{Cor}

It is a minor generalization of Theorem 2 from \cite{F}, where we got a positive solution of another Lagarias' problem
(Problem 4.4 \cite{La}).

\end{document}